\documentclass[preprint,12pt,compress]{elsarticle}

\usepackage{amssymb}
\usepackage{amsmath}
\usepackage[a4paper, body={15.5cm,23.5cm}]{geometry}
\usepackage{bbm}

\begin{document}
\journal{arXiv}
\parindent 15pt
\parskip 4pt

 \newcommand{\eps}{\varepsilon}
 \newcommand{\lam}{\lambda}
 \newcommand{\To}{\rightarrow}
 \newcommand{\as}{{\rm d}P\times {\rm d}t-a.e.}
 \newcommand{\ps}{{\rm d}P-a.s.}
 \newcommand{\jf}{\int_t^T}
 \newcommand{\tim}{\times}

 \newcommand{\F}{\mathcal{F}}
 \newcommand{\E}{\mathbf{E}}
 \newcommand{\N}{\mathbf{N}}
 \newcommand{\s}{\mathcal{S}}
 \newcommand{\T}{[0,T]}
 \newcommand{\LT}{L^2(\Omega,\F_T, P)}
 \newcommand{\Lt}{L^2(\Omega,\F_t, P)}
 \newcommand{\Ls}{L^2(\Omega,\F_s, P)}
 \newcommand{\Lr}{L^2(\Omega,\F_r, P)}
 \newcommand{\Lp}{L^p(\R^k)}
 \newcommand{\Lps}{L^2(\R^k)}
 \newcommand{\R}{{\bf R}}
 \newcommand{\RE}{\forall}

\newcommand{\gET}[1]{\underline{{\mathcal {E}}_{t,T}^g}[#1]}
\newcommand{\gEt}[1]{\underline{{\mathcal {E}}_{s,t}^g}[#1]}
\newcommand{\gETo}[1]{\underline{{\mathcal {E}}_{0,T}^g}[#1]}

\newcommand{\sgET}[1]{\overline{{\mathcal {E}}_{t,T}^g}[#1]}
\newcommand{\sgEt}[1]{\overline{{\mathcal {E}}_{s,t}^g}[#1]}
\newcommand{\sgETo}[1]{\overline{{\mathcal {E}}_{0,T}^g}[#1]}

\newcommand {\Lim}{\lim\limits_{n\rightarrow \infty}}
\newcommand {\Limk}{\lim\limits_{k\rightarrow \infty}}
\newcommand {\Limm}{\lim\limits_{m\rightarrow \infty}}
\newcommand {\llim}{\lim\limits_{\overline{n\rightarrow \infty}}}
\newcommand {\slim}{\overline{\lim\limits_{n\rightarrow \infty}}}
\newcommand {\Dis}{\displaystyle}

\def\REF#1{\par\hangindent\parindent\indent\llap{#1\enspace}\ignorespaces}

\begin{frontmatter}
\title {$L^p$ solutions of BSDEs with a new kind of
non-Lipschitz coefficients\tnoteref{fund}}
\tnotetext[fund]{Supported by the National Natural Science Foundation of China (No. 11101422) and the Fundamental Research Funds for the Central Universities (No. 2012QNA36).\vspace{0.1cm}}

\author{ShengJun FAN\corref{cor1}$^{1,2}$}
\cortext[cor1]{Corresponding author. E-mail:\ f$\_$s$\_$j@126.com.}

\author{Long JIANG$^{2}$\corref{cor}\vspace{0.2cm}}

\address{${}^1$School of Mathematical Sciences, Fudan University, Shanghai 200433, P.R. China\vspace{-0.2cm}}

\address{${}^2$College of Sciences, China University of Mining and Technology, Jiangsu 221116, P.R. China\vspace{-0.8cm}}


\begin{abstract}
In this paper, we are interested in solving multidimensional
backward stochastic differential equations (BSDEs) with a new kind
of non-Lipschitz coefficients. We establish an existence and
uniqueness result of solutions in $L^p\ (p>1)$, which
includes some known results as its particular cases.
\end{abstract}\vspace{-0.5cm}

\begin{keyword}
Backward stochastic differential equation \sep Non-Lipschitz
coefficients \sep Mao's condition\sep Constantin's condition\sep
 $L^p$ solution\vspace{0.2cm}

\MSC[2010] 60H10
\end{keyword}
\end{frontmatter}\vspace{-0.5cm}


\section {Introduction}

In this paper, we consider the following multidimensional backward
stochastic differential equation (BSDE for short in the
remaining):
\begin{equation}
    y_t=\xi+\int_t^Tg(s,y_s,z_s){\rm d}s-\int_t^Tz_s {\rm d}B_s,\ \
    t\in\T,
\end{equation}
where $T\geq 0$ is a constant called the time horizon, $\xi$ is a $k$-dimensional random vector called the terminal condition, the
random function $g(\omega,t,y,z):\Omega\tim \T\tim {\R}^{k }\tim
{\R}^{k\times d}\To {\R}^k$ is progressively measurable for each
$(y,z)$, called the generator of BSDE (1), and $B$ is a
$d$-dimensional Brownian motion. The solution
($y_{\cdot},z_{\cdot}$) is a pair of adapted processes. The triple
$(\xi,T,g)$ is called the coefficients (parameters) of BSDE (1).

Such equations, in the nonlinear case, were firstly introduced by
[5], who established an existence and uniqueness result
for solutions in $L^2$ to BSDEs under the Lipschitz assumption of the generator $g$. Since then, BSDEs have been studied with great interest, and they have gradually become an import mathematical tool in many fields such as financial mathematics, stochastic games and
optimal control, etc. In particular, many efforts have been done in
relaxing the Lipschitz hypothesis on $g$, for instance,
[3] proved the existence of a solution in $L^2$ for
(1) when $k=1$ and g is only continuous and of linear
growth in $y$ and $z$, [4] obtained an existence and uniqueness result of a solution in $L^2$ for (1) where $g$ satisfies some kind of non-lipschitz conditions, and [6] established an existence and uniqueness result of a solution in $L^2$ for (1) where $g$ satisfies some kind of monotonicity conditions in $y$. Furthermore,
[1] investigated the existence and uniqueness of a
solution in $L^p\ (p>1)$ for (1) where the generator $g$
satisfies the monotonicity condition put forward in
[6].

This paper is interested in solving multidimensional BSDEs with a
new kind of non-Lipschitz coefficients. We establish an existence
and uniqueness result of solutions in $L^p\ (p>1)$ for BSDE (1) (see Theorem 1 in Section 3), which includes the
corresponding results in [4],[2] and
[5] as its particular cases. This paper is organized as
follows. We introduce some preliminaries and lemmas in Section 2
and put forward and prove our main result in Section 3. Finally,
Section 4 is devoted to the analysis of the new kind of
non-Lipschitz coefficients, and some corollaries, remarks and
examples are also given in this section.

\section{ Preliminaries and Lemmas}

Let us first introduce some notations. First of all, let us fix
two real numbers $T\geq 0$ and $p>1$, and two positive integers $k$ and $d$. Let $(\Omega,\F,P)$ be a probability space carrying a
standard $d$-dimensional Brownian motion $(B_t)_{t\geq 0}$. Let
$(\F_t)_{t\geq 0}$ be the natural $\sigma$-algebra generated by
$(B_t)_{t\geq 0}$ and $\F=\F_T$. In this paper, the Euclidean norm
of a vector $y\in\R^k$ will be defined by $|y|$, and for an
$k\times d$ matrix $z$, we define $|z|=\sqrt{{\rm Tr}zz^*}$, where
$z^*$ is the transpose of $z$. Let $\langle x,y\rangle$ represent
the inner product of $x,y\in\R^k$. We denote by $\Lp$ the set of
all $\R^k$-valued and $\F_T$-measurable random vectors $\xi$ such
that $\E[|\xi|^p]<+\infty$, let ${\s}^p(0,T;\R^k)$ denote the set
of $\R^k$-valued, adapted and continuous processes
$(Y_t)_{t\in\T}$ such that
$$\|Y\|_{{\s}^p}:=\left( \E[\sup_{t\in\T} |Y_t|^p] \right)
^{1/p}<+\infty.$$ Moreover, let ${\rm
M}^p(0,T;\R^{k})$ (resp. ${\rm M}^p(0,T;\R^{k\times d})$) denote the
set of (equivalent classes of) $(\F_t)$-progressively measurable
${\R}^{k}$-valued (${\R}^{k\times d}$-valued) processes $(Z_t)_{
t\in\T}$ such that
$$\|Z\|_{{\rm M}^p}:=\left\{ \E\left[\left(\int_0^T |Z_t|^2\
{\rm d}t\right)^{p/2}\right] \right\}^{1/p}<+\infty.
\vspace{0.1cm}$$ Obviously, both ${\s}^p$ and ${\rm M}^p$ are Banach
spaces. As mentioned in the introduction, we will deal only with
BSDEs which are equations of type (1), where the terminal condition
$\xi$ belongs to the space $L^p(\R^k)$, and the generator $g$ is
$(\F_t)$-progressively measurable for each $(y,z)$.\vspace{0.1cm}

{\bf Definition 1}\ \ A pair of processes $(y_t,z_t)_{t\in\T}$ is
called a solution in $L^p$ to BSDE (1), if
$(y_t,z_t)_{t\in\T}\in {\s}^p(0,T;\R^{k})\times {\rm
M}^p(0,T;\R^{k\times d})$ and satisfies (1).\vspace{0.1cm}

The following Lemma 1 comes from Corollary 2.3 in [1],
which is the starting point of this paper.

{\bf Lemma 1}\ If $(y_t,z_t)_{t\in\T}$ be a solution in $L^p$
of BSDE (1), $c(p)=p[(p-1)\wedge 1]/2$ and $0\leq t\leq T$,
then
$$\begin{array}{lll}
\Dis |y_t|^p+c(p)\int_t^T |y_s|^{p-2}1_{|y_s|\neq 0}|z_s|^2\ {\rm
d}s &\leq & \Dis|\xi|^p +p\int_t^T |y_s|^{p-2}1_{|y_s|\neq
0}\langle
y_s,g(s,y_s,z_s)\rangle\ {\rm d}s\\
&& \Dis \ -p\int_t^T |y_s|^{p-2}1_{|y_s|\neq 0}\langle y_s,z_s{\rm
d}B_s\rangle. \end{array}$$

Now, let us introduce the following Proposition 1 and Proposition
2, which will play an important role in the proof of our main
result. Before that, let us first introduce the following
assumption on the generator $g$:\vspace{0.1cm}

\hspace*{-0.5cm}(A)\ \ $\as,\RE\ (y,z)\in \R^k\times\R^{k\times
d}$, $|g(\omega,t,y,z)|\leq
\psi^{{1\over p}}(|y|^p)+\lambda|z|+\varphi_t+f_t,\ \ $\vspace{0.3cm}\\
where $\lambda\geq 0$, both $\varphi_t$ and $f_t$ are two
nonnegative, $(\F_t)$-progressively measurable processes with $\E\left[\int_0^T
\varphi_t^p\ {\rm d}t\right]<+\infty$ and $\E\left[\left(\int_0^T
f_t\ {\rm d}t\right)^p\right]<+\infty$, and $\psi(\cdot):\R^+\mapsto\R^+$ is a nondecreasing and concave function with $\psi(0)=0$.\vspace{0.1cm}

{\bf Proposition 1} Let assumption (A) hold and let
$(y_t,z_t)_{t\in\T}$ be a solution in $L^p$ to BSDE (1).
Then there exists a constant $C_{\lambda,p,T}$ depending on
$\lambda$, $p$ and $T$ such that for each $t\in \T$,
$$\begin{array}{lll}
\Dis \E\left[\left(\int_t^T |z_s|^2\ {\rm
d}s\right)^{p/2}\right]&\leq & C_{\lambda,p,T}
\left\{\E\left[\sup\limits_{s\in
[t,T]}|y_s|^p\right]+\psi\left(\E\left[\sup\limits_{s\in
[t,T]}|y_s|^p\right]\right)\right.\\
&&\Dis \hspace{1.5cm}+\left.\E\left[\int_t^T \varphi_s^p\ {\rm
d}s\right]+\E\left[\left(\int_t^T f_s\ {\rm
d}s\right)^p\right]\right\}.
\end{array}$$

{\bf Proof.}\ Applying It\^{o}'s formula to $|y_t|^2$ leads to that
$$
|y_t|^2+\int_t^T |z_s|^2\ {\rm d}s=|\xi|^2+2\int_t^T \langle
y_s,g(s,y_s,z_s)\rangle \ {\rm d}s-2\int_t^T\langle y_s,z_s{\rm
d}B_s\rangle. $$ By assumption (A) we have, for
each $s\in [t,T]$,
$$
\begin{array}{lll}
2\langle y_s,g(s,y_s,z_s)\rangle & \leq & 2|y_s|\left(\psi^{{1\over
p}}(|y_s|^p)+
\lambda |z_s|+\varphi_s+f_s\right)\\
& \leq & 2\left(\sup\limits_{s\in
[t,T]}|y_s|\right)\left(\psi^{{1\over
p}}(|y_s|^p)+\varphi_s+f_s\right)+2\lambda^2 \sup\limits_{s\in
[t,T]}|y_s|^2+{|z_s|^2\over 2}.
\end{array}
$$
Thus, in view of the inequality that $2ab\leq
a^2+b^2$ we can get that
$$
\begin{array}{lll}
\Dis {1\over 2}\int_t^T |z_s|^2\ {\rm d}s &\leq & \Dis
(3+2\lambda^2T)\sup\limits_{s\in [t,T]}|y_s|^2 +\left[\int_t^T
\psi^{{1\over p}}(|y_s|^p)\ {\rm
d}s\right]^2\\
&& \Dis \ +\left[\int_t^T (\varphi_s+f_s)\ {\rm
d}s\right]^2+2\left|\int_t^T\langle y_s,z_s{\rm
d}B_s\rangle\right|.
\end{array}
$$
Then noticing that $\psi(\cdot)$ is a nondecreasing function, by
the inequality $(a+b)^{p/2}\leq 2^p(a^{p/2}+b^{p/2})$ we have
\begin{equation}
\begin{array}{lll}
\Dis \left[\int_t^T |z_s|^2\ {\rm d}s \right]^{p/2}&\leq & \Dis
c_{\lambda,p,T} \left[\sup\limits_{s\in [t,T]}|y_s|^p+
\psi(\sup\limits_{s\in [t,T]}|y_s|^p)\right.\\&&
\Dis\hspace{1.4cm} +\left.\left[\int_t^T (\varphi_s+f_s)\ {\rm
d}s\right]^p+\left|\int_t^T\langle y_s,z_s{\rm
d}B_s\rangle\right|^{p/2}\right],
\end{array}
\end{equation}
where $c_{\lambda,p,T}=2^{p+4}(3+2\lambda^2T+T^p)$ and we have used
the fact that
$$\left[\int_t^T \psi^{{1\over p}}(|y_s|^p)\
{\rm d}s\right]^p\leq T^p\psi(\sup\limits_{s\in [t,T]}|y_s|^p).$$
But by the Burkholder-Davis-Gundy (BDG) inequality, we get that for each $t\in\T$,
$$
\begin{array}{lll}
\Dis c_{\lambda,p,T}\E\left[\left|\int_t^T\langle y_s,z_s{\rm
d}B_s\rangle\right| ^{p/2} \right]&\leq & \Dis
d_{\lambda,p,T}\E\left[\left(\int_t^T |y_s|^2|z_s|^2\ {\rm
d}s\right) ^{p/4}
\right]\\
 &\leq & \Dis d_{\lambda,p,T}\E\left[\sup\limits_{s\in
[t,T]}|y_s|^{p/2}\cdot\left( \int_t^T|z_s|^2\ {\rm
d}s\right)^{p/4}\right]
\end{array}
$$ and thus,
$$
\Dis c_{\lambda,p,T}\E\left[\left|\int_t^T\langle
y_s,z_s{\rm d}B_s\rangle\right| ^{p/2} \right]\leq \Dis
{d_{\lambda,p,T}^2\over 2}\E\left[\sup\limits_{s\in
[t,T]}|y_s|^p\right]+{1\over 2}\E\left[ \left(\int_t^T|z_s|^2\
{\rm d}s\right)^{p/2}\right].
$$
Coming back to estimate (2) we get, for each $t\in \T$,
$$\begin{array}{lll}
\Dis \E\left[\left(\int_t^T |z_s|^2\ {\rm
d}s\right)^{p/2}\right]&\leq &\Dis  \bar C_{\lambda,p,T}
\left\{\E\left[\sup\limits_{s\in
[t,T]}|y_s|^p\right]+\E\left[\psi\left(\sup\limits_{s\in
[t,T]}|y_s|^p\right)\right]\right.\\
&& \Dis \hspace{1.5cm}+\left.\E\left[\left(\int_t^T
(\varphi_s+f_s)\ {\rm d}s\right)^p\right]\right\}.
\end{array}$$ Thus noticing that $\psi(\cdot)$ is a
concave function, we can deduce the desired conclusion from Jensen's inequality and H\"{o}lder's inequality. The proof is complete.\hfill $\Box$

{\bf Proposition 2} Let assumption (A) hold and let
$(y_t,z_t)_{t\in\T}$ be a solution in $L^p$ to BSDE (1).
Then there exists constants $m_p>0$ (depending on $p$) and
$K_{\lambda,p}>0$ (depending on $\lambda$ and $p$) such that for each $t\in \T$,
$$
\begin{array}{lll}
\Dis \E\left[\sup\limits_{s\in [t,T]}|y_s|^{p}\right]
&\leq &\Dis
e^{K_{\lambda,p}(T-t)}\left\{m_p\E[|\xi|^p]+m_p
\E\left[\left(\int_t^T
f_s\ {\rm
d}s\right)^{p}\right]\right.\\
&& \Dis \hspace{2.1cm}+\left.{1\over 2}\E\left[\int_t^T
\varphi_s^p\ {\rm d}s\right]+{1\over 2}\int_t^T \psi(\E[|y_s|^p])\
{\rm d}s \right\}.
\end{array}$$

{\bf Proof.}\ From Lemma 1, we get the following inequality:
$$\begin{array}{lll}
\Dis |y_t|^p+c(p)\int_t^T |y_s|^{p-2}1_{|y_s|\neq 0}|z_s|^2\ {\rm
d}s &\leq & \Dis|\xi|^p +p\int_t^T |y_s|^{p-2}1_{|y_s|\neq
0}\langle
y_s,g(s,y_s,z_s)\rangle\ {\rm d}s\\
&& \Dis \ -p\int_t^T |y_s|^{p-2}1_{|y_s|\neq 0}\langle y_s,z_s{\rm
d}B_s\rangle. \end{array}$$ Assumption (A)
yields the inequality
$$\langle
y_s,g(s,y_s,z_s)\rangle\leq |y_s|[\psi^{{1\over
p}}(|y_s|^p)+\lambda|z_s|+\varphi_s+f_s],$$ from which we deduce
that, with probability one, for each $t\in \T$,
$$\begin{array}{lll}
|y_t|^p+c(p)\int_t^T |y_s|^{p-2}1_{|y_s|\neq 0}|z_s|^2\
{\rm d}s &\leq &\Dis \Dis|\xi|^p -p\int_t^T |y_s|^{p-2}1_{|y_s|\neq
0}\langle y_s,z_s{\rm d}B_s\rangle\\
 && \Dis +p\int_t^T
|y_s|^{p-1}[\psi^{{1\over
p}}(|y_s|^p)+\lambda|z_s|+\varphi_s+f_s]\ {\rm d}s.
\end{array}$$
First of all, in view of the fact that $\psi(\cdot)$ increases at
most linearly since it is a nondecreasing concave function and
$\psi(0)=0$, we deduce from the previous inequality that, $\ps$,
$$\int_0^T |y_s|^{p-2}1_{|y_s|\neq 0}|z_s|^2\
{\rm d}s<+\infty.$$ Moreover, making use of Young's inequality (
$a^rb^{1-r}\leq ra+(1-r)b$ for each $a\geq 0$, $b\geq 0$ and $r\in
(0,1)$) and the inequality $(a+b)^p\leq 2^p(a^p+b^p)$ we can
obtain
$$\begin{array}{lll}
\hspace*{-0.5cm}p\int_t^T |y_s|^{p-1}(\psi^{{1\over
p}}(|y_s|^p)+\varphi_s)\ {\rm d}s &= & \Dis  p\int_t^T \left\{\left(\theta^{1\over
p-1}|y_s|^{p}\right)^{p-1\over p}\cdot \left[{1\over
\theta}\left(\psi^{{1\over p}}(|y_s|^p)+\varphi_s\right)^{p}
\right]^{1\over
p}\right\}\ {\rm d}s\\
&\leq & \Dis  (p-1)\theta^{1\over p-1}\int_t^T |y_s|^{p}\ {\rm
d}s+{2^p\over \theta} \int_t^T (\psi(|y_s|^p)+\varphi_s^p)\ {\rm
d}s,
\end{array}$$
where $\theta>0$ will be chosen later. And, from the inequality
that $ab\leq (a^2+b^2)/2 $ we get that
$$
\begin{array}{lll}
\Dis p\lambda |y_s|^{p-1}|z_s|&=&\Dis p\left({\sqrt{2}\lambda
\over \sqrt{1\wedge (p-1)}}|y_s|^{p\over 2}\right)
\left(\sqrt{{1\wedge (p-1)}\over 2}|y_s|^{p-2\over 2}|z_s|\right)\\
&\leq & \Dis {p\lambda^2\over 1\wedge (p-1)}|y_s|^p +{c(p)\over
2}|y_s|^{p-2}1_{|y_s|\neq 0}|z_s|^2.
\end{array}$$
Thus for each $t\in \T$, we
have
\begin{equation}
\Dis |y_t|^p+{c(p)\over 2}\int_t^T |y_s|^{p-2}1_{|y_s|\neq
0}|z_s|^2\ {\rm d}s \leq  \Dis X_t-p\int_t^T
|y_s|^{p-2}1_{|y_s|\neq 0}\langle y_s,z_s{\rm d}B_s\rangle,
\end{equation}
where $X_t=|\xi|^p +d_{\lambda,p,\theta} \int_t^T |y_s|^{p}\ {\rm
d}s+
 {2^p\over \theta} \int_t^T (\psi(|y_s|^p)+\varphi_s^p)\ {\rm d}s+
 p\int_t^T |y_s|^{p-1}f_s\ {\rm d}s
$ with $d_{\lambda,p,\theta}=(p-1)\theta^{1/(p-1)}
+{p\lambda^2/[1\wedge (p-1)]}>0$.

It follows from the BDG inequality that $\{M_t:=\int_0^t
|y_s|^{p-2}1_{|y_s|\neq 0}\langle y_s,z_s{\rm
d}B_s\rangle\}_{t\in\T}$ is a uniformly integrable martingale. In
fact, we have, by Young's inequality,
$$
\begin{array}{lll}
\Dis \E\left[\langle M, M \rangle^{1/2}_T\right]&\leq &
\Dis\E\left[\sup\limits_{s\in [0,T]}|y_s|^{p-1}\cdot\left(
\int_0^T|z_s|^2\ {\rm
d}s\right)^{1/2}\right]\\
&= & \Dis \E\left\{\left(\sup\limits_{s\in
[0,T]}|y_s|^{p}\right)^{p-1\over p}\cdot \left[\left(\int_0^T
|z_s|^2\
{\rm d}s\right)^{p/2}\right]^{1\over p}\right\}\\
&\leq &\Dis {(p-1)\over p}\E\left[\sup\limits_{s\in
[0,T]}|y_s|^p\right]+{1\over p}\E\left[\left( \int_0^T|z_s|^2\
{\rm d}s\right)^{p/2}\right]<+\infty.
\end{array}
$$
Coming back to inequality (3), and taking the expectation, we
get both
\begin{equation}
{c(p)\over 2}\E\left[\int_t^T |y_s|^{p-2}1_{|y_s|\neq 0}|z_s|^2\
{\rm d}s\right] \leq \E[X_t]
\end{equation}
and
\begin{equation}
\E\left[\sup\limits_{s\in [t,T]}|y_s|^{p}\right]\leq
\E[X_t]+k_p\E\left[\left(\langle M, M \rangle_T-\langle M, M
\rangle_t\right) ^{1/2}\right],
\end{equation}\vspace{0.1cm}
where we have used the BDG inequality in the last inequality.

On the other hand, making use of Young's inequality we have also
$$
\begin{array}{lll}
&&\Dis k_p\E\left[\left(\langle M, M \rangle_T-\langle M, M
\rangle_t\right) ^{1/2}\right]\\
&\leq & \Dis k_p\E\left[\sup\limits_{s\in
[t,T]}|y_s|^{p/2}\cdot\left( \int_t^T|y_s|^{p-2}1_{|y_s|\neq
0}|z_s|^2\ {\rm
d}s\right)^{1/2}\right]\\
&\leq &\Dis {1\over 2}\E\left[\sup\limits_{s\in
[t,T]}|y_s|^p\right]+{k_p^2\over 2}\E\left[
\int_t^T|y_s|^{p-2}1_{|y_s|\neq 0}|z_s|^2\ {\rm d}s\right].
\end{array}
$$
Coming back to inequalities (4) and (5), we obtain
$$\E\left[\sup\limits_{s\in
[t,T]}|y_s|^{p}\right]\leq k'_p\E[X_t].$$ Applying once again
Young's inequality, we get
$$\begin{array}{lll}
\Dis pk'_p\E\left[\int_t^T |y_s|^{p-1}f_s\ {\rm d}s \right]&\leq
&\Dis pk'_p\E\left[\sup\limits_{s\in
[t,T]}|y_s|^{p-1}\int_t^T f_s\ {\rm d}s \right]\\
&=& \Dis \E\left\{\left({p\over 2(p-1)} \sup\limits_{s\in
[t,T]}|y_s|^{p}\right)^{p-1\over p}\cdot \left[{pk''_p\over
2}\left(\int_t^T f_s\ {\rm
d}s\right)^{p}\right]^{1\over p}\right\}\\
&\leq & \Dis  {1\over 2}\E\left[\sup\limits_{s\in
[t,T]}|y_s|^{p}\right]+{k''_p\over 2}\E\left[\left(\int_t^T f_s\
{\rm d}s\right)^{p}\right],
\end{array}$$
from which we deduce, combing back to the definition of $X_t$,
that
$$\begin{array}{lll}
\Dis \E\left[\sup\limits_{s\in [t,T]}|y_s|^{p}\right] &\leq &
\Dis
2k'_p\E\left[|\xi|^p +d_{\lambda,p,\theta} \int_t^T |y_s|^{p}\
{\rm d}s+
 {2^p\over \theta}\int_t^T (\psi(|y_s|^p)+\varphi_s^p)
 \ {\rm d}s\right]\\
 && \Dis\  +k''_p\E\left[\left(\int_t^T f_s\
{\rm d}s\right)^{p}\right].
\end{array}$$
By letting $\theta=2^{p+2}k'_p$ and $h_t=\E\left[\sup\limits_{s\in
[t,T]}|y_s|^{p}\right]$ in the previous inequality and using Fubini theorem and Jensen's inequality, noticing that $\psi(\cdot)$ is a
concave function, we have, for each $t\in \T$,
$$\begin{array}{lll}
\Dis h_t &\leq & \Dis  2k'_p\E[|\xi|^p]
 +k''_p\E\left[\left(\int_t^T f_s\
{\rm d}s\right)^{p}\right]+{1\over 2}\E\left[\int_t^T \varphi_s^p\
{\rm d}s\right]\\
&& \Dis +{1\over 2}\int_t^T \psi(\E[|y_s|^p])\ {\rm
d}s+2k'_pd_{\lambda,p,\theta}\int_t^T h_s\ {\rm d}s.
\end{array}
$$
Finally, Gronwall's inequality yields that for each $t\in \T$,
$$\begin{array}{lll}
\Dis h_t&\leq & \Dis
e^{2k'_pd_{\lambda,p,\theta}(T-t)}\left\{2k'_p\E[|\xi|^p]
 +k''_p\E\left[\left(\int_t^T f_s\
{\rm d}s\right)^{p}\right]\right.\\
&& \Dis \hspace{2.7cm}+\left.{1\over 2}\E\left[\int_t^T
\varphi_s^p\ {\rm d}s\right]+
 {1\over 2}\int_t^T \psi(\E[|y_s|^p])\ {\rm d}s\right\}.
 \end{array}
 $$
Then we complete the proof of Proposition 2.\hfill $\Box$

\section{ Main Result and Its Proof }

In this section, we will put forward and prove our main result.
Let us first introduce the following assumptions:\vspace{0.1cm}

(H1) \ There exists a nondecreasing and concave function
$\rho(\cdot):\R^+\mapsto \R^+$ with $\rho(0)=0$, $\rho(u)>0$ for
$u>0$ and $\int_{0^+} {{\rm d}u\over \rho(u)}=+\infty$ such that
$\as,$
$$\RE y_1,y_2\in \R^k,z\in\R^{k\times d},\ \
|g(\omega,t,y_1,z)-g(\omega,t,y_2,z)|^p\leq \rho(|y_1-y_2|^p).$$

(H2)\ There exists a constant $C\geq 0$ such that $\as,$
$$\RE y\in \R^k,z_1,z_2\in\R^{k\times d},\ \
|g(\omega,t,y,z_1)-g(\omega,t,y,z_2)|\leq
C|z_1-z_2|.$$

(H3)\  $\Dis\E\left[\left(\int_0^T |g(t,0,0)|\ {\rm
d}t\right)^p\right]<+\infty$.\vspace{0.2cm}

{\bf Remark 1}\ Since $\rho(\cdot)$ is a nondecreasing and concave
function with $\rho(0)=0$, it increases at most linearly, i.e.,
there exists a constant $A>0$ such that $\rho(x)\leq A(x+1)$ for
each $x\geq 0$.\vspace{0.2cm}

The following Theorem 1 is the main result of this
paper.

{\bf Theorem 1}\ \  Let $g$ satisfy assumptions (H1)-(H3). Then
for each $\xi\in\Lp$, the BSDE with parameters $(\xi,T,g)$ has a
unique solution in $L^p$.

We can construct the Picard approximate sequence of the BSDE with
parameters $(\xi,T,g)$ as follows:
\begin{equation}
y_t^0=0;\ \ \ \  \Dis y_t^n=\xi+\int_t^Tg(s,y_s^{n-1},z_s^n){\rm
d}s-\int_t^Tz_s^n {\rm d}B_s,\ \  t\in\T.
\end{equation}
Indeed, for each $n\geq 1$, by (H1) and Remark 1 we have
$$\begin{array}{lll}
|g(s,y_s^{n-1},0)|&\leq & \Dis |g(s,0,0)|+\rho^{1\over p}
(|y_s^{n-1}|^p)\\
&\leq & \Dis |g(s,0,0)|+A^{1\over p}(|y_s^{n-1}|+1),
\end{array}$$ and then
$$
\begin{array}{lll}
\Dis \E\left[\left(\int_0^T|g(s,y_s^{n-1},0)|\ {\rm
d}s\right)^p\right] &\leq & \Dis
2^p\E\left[\left(\int_0^T|g(s,0,0)|\ {\rm
d}s\right)^p\right]\\
&& \Dis +A(2T)^p\left(\E\left[\sup\limits_{s\in\T}
|y_s^{n-1}|^p\right] +1\right).
\end{array}
$$
Then the generator $g(t,y_t^{n-1},z)$ of BSDE (6) satisfies (H2)
and (H3). It follows from Theorem 4.2 in [1] that, for
each $n\geq 1$, equation (6) has a unique solution
$(y_t^n,z_t^n)_{t\in\T}$ in $L^p$. Concerning the
processes $(y_t^n,z_t^n)_{t\in\T}$, we have the following Lemma 2
and Lemma 3.

{\bf Lemma 2}\ Under the hypotheses of Theorem 1, there exists a
constant $c_1>0$ depending only on $C$ and $p$, and a constant
$K>0$ depending only on $C$, $p$ and $T$, such that for each $t\in
\T, n,m\geq 1$,
\begin{equation}
\Dis \E\left[\sup\limits_{s\in [t,T]}|y_s^{n+m}-y_s^n|^p\right]
\leq \Dis {1\over 2}e^{c_1(T-t)}\int_{t}^T\rho\left(\E\left[
|y_s^{n+m-1}-y_s^{n-1}|^p\right]\right)\ {\rm d}s.
\end{equation}
and
\begin{equation}
\begin{array}{lll}
\Dis \E\left[\left(\int_t^T |z_s^{n+m}-z_s^n|^2\ {\rm
d}s\right)^{p/2}\right]&\leq&\Dis K
\left\{\E\left[\sup\limits_{s\in
[t,T]}|y_s^{n+m}-y_s^n|^p\right]\right. \\
&& \Dis \hspace{0.8cm}+\left.\int_t^T \rho\left(\E\left[
|y_s^{n+m-1}-y_s^{n-1}|^p\right]\right)\ {\rm d}s\right\}.
\end{array}
\end{equation}

{\bf Proof.} It follows from (6) that the process
$(y_t^{n+m}-y_t^n,z_t^{n+m}-z_t^n)_{t\in\T}$ is a solution in
$L^p$ of the following BSDE
\begin{equation}
\Dis y_t=\int_t^T f_{n,m}(s,z_s){\rm d}s-\int_t^Tz_s {\rm d}B_s,\
\ t\in\T
\end{equation}
where
$$f_{n,m}(s,z):= g(s,y_s^{n+m-1},z+z_s^{n})-g(s,y_s^{n-1},
z_s^{n}).$$
By (H1) and (H2) we have
$$|f_{n,m}(s,z)|\leq \rho^{1\over p}(|y_s^{n+m-1}-y_s^{n-1}|^p)
+C|z|,$$ which means that assumption (A) is satisfied for the
generator $f_{n,m}(t,z)$ of BSDE (9) with  $\psi(\cdot)\equiv
0$, $\lambda=C$, $f_t\equiv 0$ and $\varphi_t=\rho^{1\over
p}(|y_t^{n+m-1}-y_t^{n-1}|^p)$ by Remark 1. Thus, the conclusions
(7) and (8) follows, in view of the fact that $\rho(\cdot)$ is a
concave function, from Proposition 2 and Proposition 1, and then
Fubini theorem and Jensen's inequality.\hfill$\Box$\vspace{0.1cm}

{\bf Lemma 3}\ Under the hypotheses of Theorem 1, there exists
$T_1\in [0,T)$ independent of the terminal condition $\xi$ and a
constant $M\geq 0$ such that for each $n\geq 1$, $t\in [T_1,T]$,
$$\E\left[\sup_{r\in
[t,T]}|y_r^n|^p\right]\leq M.$$

{\bf Proof.}\ Making use of the hypotheses of Theorem 1 we know
that
$$
\begin{array}{lll}
|g(s,y_s^{n-1},z)| &\leq & |g(s,y_s^{n-1},z)-g(s,0,0)|+
|g(s,0,0)|\\
 &\leq &
\rho^{1\over p}(|y_s^{n-1}|^p)+C|z|+|g(s,0,0)|
\end{array}$$
Thus, assumption (A) is satisfied for the generator
$g(t,y_t^{n-1},z)$ of BSDE (6) with $\psi(\cdot)\equiv 0$,
$\lambda=C$, $f_t=|g(t,0,0)|$ by (H3) and $\varphi_t=\rho^{1\over
p}(|y_t^{n-1}|^p)$ by Remark 1. Consequently, in view of the fact
that $\rho(\cdot)$ is a concave function, from Proposition 2 and
then Fubini theorem and Jensen's inequality we get that there
exist two positive constants $c_2$ and $c_3$ depending only on
$C$ and $p$ such that for each $n\geq 1$ and each $t\in [0,T]$,
\begin{equation}
\E\left[\sup\limits_{r\in [t,T]}|y_r^n|^p\right]\leq  \Dis
\mu_t+{1\over 2}e^{c_3(T-t)}\int_{t}^T
\rho\left(\E\left[|y_s^{n-1}|^p\right]\right) \ {\rm d}s
\end{equation}
where
$\mu_t=c_2e^{c_3(T-t)}\left\{\E|\xi|^p+\E\left[\left(\int_{t}^T|g(s,0,0)|{\rm
d}s\right)^p\right]\right\}\geq 0.$

Now, let $M=2\mu_0+2AT$ and $T_1=\max\{T-{\ln 2/c_1},T-{\ln
2/c_3},T-{1/ 2A},0\}$ where $c_1$ is defined in Lemma 2 and $A$ is
defined in Remark 1. Then for each $t\in [T_1,T]$, we have
\begin{equation}
{1\over 2}e^{c_1(T-t)}\leq 1,\ {1\over 2}e^{c_3(T-t)}\leq 1 \ \
{\rm and}\  \ A(T-t)\leq {1\over 2}.
\end{equation}
 Thus, from (10) and (11) we have
\begin{equation}
\E\left[\sup\limits_{r\in [t,T]}|y_r^n|^p\right]\leq  \Dis
\mu_0+\int_{t}^T \rho\left(\E\left[|y_s^{n-1}|^p\right]\right) \
{\rm d}s,\ \ t\in [T_1,T].
\end{equation}
Since $\rho(\cdot)$ is a nondecreasing function, by (12), Remark 1
and (11) we can deduce that, for each $t\in [T_1,T]$,
$$
\Dis\E\left[\sup\limits_{r\in [t,T]}|y_r^1|^p\right]
\leq \mu_0\leq
M,
$$
$$
\Dis\E\left[\sup\limits_{r\in [t,T]}|y_r^2|^p\right]\leq
\mu_0+\int_{t}^T\rho(M)\ {\rm d}s \leq \mu_0+A(M+1)(T-t) \leq
\mu_0+{M\over 2}+AT=M,
$$
$$
\Dis\E\left[\sup\limits_{r\in [t,T]}|y_r^3|^p\right]\leq
\mu_0+\int_{t}^T\rho(M)\ {\rm d}s\leq \mu_0+A(M+1)(T-t) \leq
\mu_0+{M\over 2}+AT=M.
$$
By induction, we can prove that for all $n\geq 1$ and all $t\in
[T_1,T]$,
$$\E\left[\sup_{r\in
[t,T]}|y_r^n|^p\right]\leq M.$$ The proof of Lemma 3 is complete.\hfill $\Box$\vspace{0.1cm}

With the help of Lemma 2 and Lemma 3, we can prove Theorem
1.

{\bf Proof of Theorem 1.}\ Existence:\ Define a sequence of
functions $\{\varphi_n(t)\}_{n\geq 1}$ as follows:
\begin{equation}
\Dis \varphi_0(t)=\int_t^T\rho(M)\ {\rm d}s; \ \ \Dis
\varphi_{n+1}(t)=\int_t^T\rho(\varphi_n(s))\ {\rm d}s.
\end{equation}
Then for all $t\in [T_1,T]$, from the proof of Lemma 3 we have
$$
\begin{array}{l}
\Dis \varphi_0(t)=\int_{t}^T\rho(M)\ {\rm d}s\leq M,\\
\Dis \varphi_1(t)=\int_{t}^T \rho(\varphi_0(s))\ {\rm d}s
\leq \int_{t}^T\rho(M)\ {\rm d}s=\varphi_0(t)\leq M,\\
\Dis \varphi_2(t)=\int_{t}^T \rho(\varphi_1(s))\ {\rm d}s \leq
\int_{t}^T\rho(\varphi_0(s))\ {\rm d}s=\varphi_1(t)\leq M.
\end{array}
$$
By induction, we can prove that for all $n\geq 1$, $
\varphi_n(t)$ satisfies
$$0\leq \varphi_{n+1}(t)\leq \varphi_n(t)\leq \cdots\leq
\varphi_{1}(t)\leq \varphi_0(t)\leq M.$$ Then, for each $t\in
[T_1,T]$, the limit of the sequence $\{\varphi_n(t)\}_{n\geq 1}$
must exist, we denote it by $\varphi(t)$. Thus, letting $n\To
\infty$ in (13), in view of the facts that $\rho(\cdot)$ is a
continuous function and $\rho(\varphi_n(s))\leq \rho(M)$ for each
$n\geq 1$, we can deduce from the Lebesgue dominated convergence
theorem that for each $t\in [T_1,T]$,
$$\varphi(t)=\int_t^T\rho(\varphi(s))\ {\rm d}s.$$
Then Bihari's inequality (see Lemma 3.6 in [4]) yields
that for each $t\in [T_1,T]$, $\varphi(t)=0$.

Now, for all $t\in [T_1,T]$, $n,m\geq 1$, thanks to Lemma 3, (7)
in Lemma 2 and inequality (11) we have,
$$
\hspace*{-1.7cm}\begin{array}{l}
\E\left[\sup\limits_{r\in [t,T]}|y_r^n|^p\right]\leq M,\\
\Dis \E\left[\sup\limits_{r\in [t,T]}|y_r^{1+m}-y_r^1|^p\right]
\leq \Dis \int_t^T \rho\left(\E\left[|y_s^m|^p\right]\right)\ {\rm
d}s
\leq \Dis \int_t^T\rho(M)\ {\rm d}s =\varphi_0(t)\leq M,\\
\end{array}
$$
$$
\begin{array}{l}
\Dis \E\left[\sup\limits_{r\in
[t,T]}|y_r^{2+m}-y_r^2|^p\right]\leq \Dis
\int_t^T\rho\left(\E\left[|y_s^{1+m}-y_s^1|^p\right]\right)\ {\rm
d}s\leq \Dis \int_t^T\rho(\varphi_0(s))\ {\rm d}s=\varphi_1(t)\leq
M,\\
\Dis \E\left[\sup\limits_{r\in
[t,T]}|y_r^{3+m}-y_r^3|^p\right]\leq \Dis \int_t^T
\rho\left(\E\left[|y_s^{2+m}-y_s^2|^p\right]\right)\ {\rm d}s \leq
\Dis \int_t^T\rho(\varphi_1(s))\ {\rm d}s=\varphi_2(t) \leq M.
\end{array}\vspace{0.2cm}
$$
By induction, we can derive that
$$\E\left[\sup\limits_{T_1\leq r\leq T}|y_r^{n+m}-y_r^n|^p\right]
\leq \varphi_{n-1}(T_1)\rightarrow 0,\ \ n\rightarrow \infty.
$$
which means that $\{y_t^n\}_{n\geq 1}$ is a Cauchy sequence in
$\s^p(T_1,T;\R^k)$. Furthermore, in view of the fact that
$\rho(\cdot)$ is a continuous function we also know from (8) in
Lemma 2 that $\{z_t^n\}_{n\geq 1}$ is a Cauchy sequence in ${\rm
M}^p(T_1,T;\R^{k\times d})$. Define their limits by $(y_t)_{t\in
[T_1,T]}$ and $(z_t)_{t\in [T_1,T]}$ respectively. Thus, by letting
$n\To \infty$ in (6), we get that $(y_t,z_t)$ is a solution in
$L^p$ to the BSDE with parameters $(\xi,T,g)$ on $[T_1,T]$.
Note from Lemma 3 that the $T_1$ does not depend on the terminal
condition $\xi$. Hence we can deduce by iteration the existence on
$[T-l(T-T_1),T]$, for each $l$, and therefore the existence on the
whole $[0,T]$. The existence has been proved.

Uniqueness: Let $(y_t^i,z_t^i)_{t\in\T}\ (i=1,2)$ be two solutions
in $L^p$ of the BSDE with parameters $(\xi,T,g)$. It follows
that $(y_t^1-y_t^2,z_t^1-z_t^2)_{t\in\T}$ is a solution in $L^p$ to the following BSDE:
\begin{equation}
\Dis y_t=\int_t^T \hat g(s,y_s,z_s){\rm d}s-\int_t^Tz_s {\rm
d}B_s,\ \ t\in\T,
\end{equation}
where
$$\hat g(s,y,z):= g(s,y+y_s^{2},z+z_s^{2})-g(s,y_s^{2},
z_s^{2}).\vspace{0.2cm}$$
By (H1) and (H2) we have
$$|\hat g(s,y,z)|\leq \rho^{1\over p}(|y|^p)
+C|z|,$$ which means that assumption (A) is satisfied for the
generator $\hat g(t,y,z)$ of BSDE (14) with
$\psi(\cdot)=\rho(\cdot)$, $\lambda=C$, $\varphi_t\equiv 0$ and
$f_t\equiv 0$. Then from Proposition 2 and Proposition 1 we can
obtain that there exists a constant $c_4>0$ depending only on $C$
and $p$, and a constant $c_5>0$ depending only on $C$, $p$ and $T$,
such that for $t\in [0,T]$,
\begin{equation}
\Dis \E\left[|y_t^{1}-y_t^2|^p\right]\leq  \Dis {1\over
2}e^{c_4(T-t)}\int_{t}^T\rho\left(\E\left[
|y_s^{1}-y_s^{2}|^p\right]\right)\ {\rm d}s
\end{equation}
and
\begin{equation}
\Dis\E\left[\left(\int_t^T |z_s^1-z_s^2|^2\ {\rm
d}s\right)^{p/2}\right]\leq  c_5 \left\{\E\left[\sup\limits_{s\in
[t,T]}|y_s^1-y_s^2|^p\right]+\rho\left(\E\left[\sup\limits_{s\in
[t,T]}|y_s^1-y_s^2|^p\right]\right)\right\}.\
\end{equation}
Then from (15) Bihari's inequality (see Lemma 3.6 in
[4]) yields that for each $t\in [0,T]$,
$\E\left[|y_t^{1}-y_t^2|^p\right]=0$. This means $y_t^1=y_t^2$ for
all $t\in [0,T]$ almost surely. We can immediately deduce that
$z_t^1=z_t^2$ for all $t\in [0,T]$ almost surely by (16). The
proof of the Theorem 1 is then complete.\hfill $\Box$

\section{ Corollaries, Remarks and Examples }

In this section, we are devoted to the analysis of the new kind of
non-Lipschitz coefficients. Some corollaries, remarks and examples
are given to show that Theorem 1 of this paper is a generalization
of the corresponding results in [4],[2] and
[5]. Firstly, by Theorem 1 the following corollary is
immediate. By H\"{o}lder's inequality we know that it generalizes
Theorem 2.1 in [4] where (H3) is replaced with
$g(\cdot,0,0)\in {\rm M}^2(0,T;\R^{k})$.\vspace{0.1cm}

{\bf Corollary 1}\ \ Let $g$ satisfy (H1) with $p=2$, (H2) and (H3). Then for each $\xi\in\Lps$, the BSDE with parameters $(\xi,T,g)$ has a unique solution in $L^2$.\vspace{0.1cm}

Furthermore, by letting $\rho(x)=\mu^px$ with $\mu>0$ in (H1) we can obtain the following classical Lipschitz assumption in $y$ with
respect to the generator $g$:\vspace{0.1cm}

(H1') \ There exists a constant $\mu\geq 0$ such that $\as,$
$$\RE y_1,y_2\in \R^k,z\in\R^{k\times d},\ \
|g(\omega,t,y_1,z)-g(\omega,t,y_2,z)|\leq \mu|y_1-y_2|.$$
Consequently, by Theorem 1 we can get the
following corollary, which generalizes the main result in
[5] where $p=2$.\vspace{0.1cm}

{\bf Corollary 2}\ \ Let $g$ satisfy (H1'), (H2) and (H3). Then for
each $\xi\in\Lp$, the BSDE with parameters $(\xi,T,g)$ has a unique
solution in $L^p$.\vspace{0.1cm}

{\bf Remark 2} In the following, we will show that the
concavity condition of $\rho(\cdot)$ in (H1) can be actually lifted
and that the bigger the $p$, the stronger the (H1). To be
precise, we need to prove that if $g$ satisfy the following
assumption (H1'') with $q\geq p$, then $g$ must satisfy
(H1).\vspace{0.1cm}

(H1'') \ There exists a nondecreasing and continuous function
$\kappa(\cdot):\R^+\mapsto \R^+$ with $\kappa(0)=0$, $\kappa(u)>0$
for $u>0$ and $\int_{0^+} {{\rm d}u\over \kappa(u)}=+\infty$ such
that $\as,$
$$\RE y_1,y_2\in \R^k,z\in\R^{k\times d},\ \
|g(\omega,t,y_1,z)-g(\omega,t,y_2,z)|^q\leq \kappa(|y_1-y_2|^q).$$

In order to show Remark 2, we need the following technical
Lemma. Its proof can be done by means of approximation procedures in [2], here we omit it.

{\bf Lemma 4} Let $\rho(\cdot)$ be a nondecreasing and concave
function on $\R^+$ with $\rho(0)=0$. Then we have
\begin{equation}
\RE\ r>1,\ \ \rho^r(x^{1\over r})\ {\rm is\ also\  a\ nondecreasing\
and \ concave\ function\ on }\ \R^+.
\end{equation}
Moreover, if $\rho(u)>0$ for $u>0$ and $\int_{0^+} {{\rm d}u\over
\rho(u)}=+\infty$, then
\begin{equation}
\RE\ r<1,\ \ \ \int_{0^+} {{\rm d}u\over \rho^r(u^{1\over
r})}=+\infty.
\end{equation}

Now, we can show that (H1'')$\Longrightarrow$(H1). Let us assume
that (H1'') holds for $g$. Then we have, $\as,$
$$
\RE y_1,y_2\in \R^k,z\in\R^{k\times d},\ \
|g(\omega,t,y_1,z)-g(\omega,t,y_2,z)|\leq
\rho_1(|y_1-y_2|),
$$ where
$\rho_1(u):=\kappa^{{1\over q}}(u^q)$. Obviously, $\rho_1(u)$ is a
continuous and nondecreasing function on $\R_+$ with $\rho_1(0)=0$
and $\rho_1(x)>0$ for $x>0$, but it is not necessary to be concave.
However, it follows from [2] that if $g$ satisfies the
above condition, then there exists a concave and nondecreasing function $\rho_2(\cdot)$ such that $\rho_2(0)=0$, $\rho_2(u)\leq 2\rho_1(u)$ for $u\geq 0$, and $\as,$
$$
\RE y_1,y_2\in \R^k,z\in\R^{k\times d},\ \
|g(\omega,t,y_1,z)-g(\omega,t,y_2,z)|\leq \rho_2(|y_1-y_2|).
$$
Thus, $\as,$
$$
\RE y_1,y_2\in \R^k,z\in\R^{k\times d},\ \
|g(\omega,t,y_1,z)-g(\omega,t,y_2,z)|^p\leq \bar\rho(|y_1-y_2|^p),
$$
where $\bar\rho(u):=\rho_2^p(u^{1\over p})+u$. It is clear that
$\bar\rho(0)=0$ and $\bar\rho(u)>0$ for $u>0$. Moreover, it
follows from (17) in Lemma 4 that $\bar\rho(\cdot)$ is also a
nondecreasing and concave function since $p>1$ and $\rho_2(\cdot)$
is a nondecreasing and concave function. Thus, to prove that (H1)
holds, it suffices to show that $\int_{0^+}{{\rm d}u\over
\bar\rho(u)}=+\infty$. Indeed, if $\rho_2(1)=0$, then since
$\rho_2(u)=0$ for each $u\in [0,1]$, we have
$$\int_{0^+}{{\rm d}u\over \bar\rho(u)}=\int_{0^+}{{\rm d}u\over u}=
+\infty.$$

On the other hand, if $\rho_2(1)>0$, since $\rho_2(\cdot)$ is a concave function with $\rho_2(0)=0$, we know that $\RE u\in [0,1],\ \rho_2(u)\geq u\rho_2(1)$, and then
$$\rho_2^{p}({u^{1\over p}})\geq \left(u^{1\over p} \rho_2(1)
\right)^{p}=\rho_2^{p}(1)u.$$
Thus, we have
\begin{equation}
\RE\ u\geq 0,\ \bar\rho(u)=\rho_2^p(u^{1\over p})+u \leq
K\rho_2^p(u^{1\over p})\leq K2^p\rho_1^p(u^{1\over p})
=K2^p\kappa^{p\over q}({u^{q\over p}}),
\end{equation}
where $K=1+{1/\rho_2^{p}(1)}.$ Consequently, if $q=p$, then
$$\int_{0^+}{{\rm d}u\over \bar\rho(u)}\geq {1\over
K 2^p} \int_{0^+}{{\rm d}u\over
 \kappa(u)} =+\infty.
$$
Thus, we have proved that (H1'') with $q=p$ implies (H1). Hence,
now we can assume that the $\bar \kappa (\cdot)$ in (H1'') is a
concave function. Then, if $q>p$, from (19) and (18) in Lemma 4
with $\rho(\cdot)=\kappa(\cdot)$ and $r={p/q}<1$ we have
$$\int_{0^+}{{\rm d}u\over \bar\rho(u)}\geq {1\over
K 2^p} \int_{0^+}{{\rm d}u\over \kappa^{p\over q}({u^{q\over p}})}
=+\infty.
$$
Thus, (H1) holds. Hence (H1'')$\Longrightarrow$(H1), i.e.,
the concavity condition of $\rho(\cdot)$ in (H1) can be actually
lifted and the bigger the $p$, the stronger the (H1).\vspace{0.1cm}

Furthermore, let us introduce a stronger assumption (H1*) than
(H1):\vspace{0.1cm}

(H1*) \ There exists a nondecreasing and continuous function
$\kappa(\cdot):\R^+\mapsto \R^+$ with $\kappa(0)=0$, $\kappa(u)>0$
for $u>0$ and $\int_{0^+} {u^{p-1}\over \kappa^p(u)}{\rm
d}u=+\infty$ such that $\as,$
$$\RE y_1,y_2\in \R^k,z\in\R^{k\times d},\ \
|g(\omega,t,y_1,z)-g(\omega,t,y_2,z)|\leq \kappa(|y_1-y_2|).$$

In the following, we show (H1*)$\Longrightarrow$ (H1). In fact, if $g$ satisfies (H1*), then we have, $\as,$
$$\RE y_1,y_2\in \R^k,z\in\R^{k\times d},\ \
|g(\omega,t,y_1,z)-g(\omega,t,y_2,z)|^p\leq \rho(|y_1-y_2|^p),$$
where $\rho(u)=\kappa^p(u^{1\over p})$. And, we have also,
\begin{equation}
\int_{0+}{{\rm d}u\over \rho(u)}=\int_{0+}{{\rm d}u \over
\kappa^p(u^{1\over p})}=\int_{0^+} {pu^{p-1}\over \kappa^p(u)}{\rm
d}u=+\infty.
\end{equation}
Thus, in view of Remark 2, we know that (H1) is true.
Therefore, from Theorem 1 the following corollary is immediate. And, by H\"{o}lder's inequality we know that it generalizes the
corresponding result in [2], where $p=2$ and (H3) is
replaced by $g(\cdot,0,0)\in {\rm
M}^2(0,T;\R^{k})$:\vspace{0.1cm}

{\bf Corollary 3}\ \ Let $g$ satisfy (H1*), (H2) and (H3). Then
for each $\xi\in\Lp$, the BSDE with parameters $(\xi,T,g)$ has a
unique solution in $L^p$.\vspace{0.1cm}

{\bf Remark 3}\ According to the classical theory of uniformly
continuous functions, we can assume that the $\kappa(\cdot)$ in
(H1*) is a concave function. Thus, applying (18) in Lemma 2, by
letting $\rho(u)=\kappa^q(u^{1/q})$ and $r=p/q$ with $q>p$ we
deduce that if $\int_{0+}{{\rm d}u \over \kappa^q(u^{1/q})}=+\infty$, then $\int_{0+}{{\rm d}u \over \kappa^p(u^{1/p})}=+\infty$. Thus, noticing (20) we know that the bigger the
$p$, the stronger the (H1*).\vspace{0.1cm}

To the end, we give an example.

{\bf Example 1} Let $g(t,y,z)=h(|y|)+|z|+|B_t|$, where
$$h(x):=x|\ln x|^{1/p}\cdot 1_{0<x\leq \delta}+(h'(\delta-)
(x-\delta)+h(\delta))\cdot 1_{x>\delta}$$
with $\delta>0$ small enough. It is clear that $g$ satisfies (H2) and (H3). We can also prove that $g$ satisfies (H1*) with $\kappa(\cdot)=h(\cdot)$ by verifying that $\int_{0^+}{u^{p-1}\over h^p(u)}{\rm d}u=+\infty$, $h(\cdot)$
is a sub-additive function and then the following inequality holds:
$$\RE y_1,y_2\in \R^k,z\in\R^{k\times d},\ \
|g(\omega,t,y_1,z)-g(\omega,t,y_2,z)|\leq
h(|y_1-y_2|).$$
Thus, this generator $g$ satisfies all conditions in Corollary 3. Consequently, the BSDE with parameters $(\xi,T,g)$ has a unique solution in $L^p$ for each $\xi\in\Lp$.

Finally, it is worth mentioning that we can directly verify that
for each $q>p$, $\int_{0^+}{u^{q-1}\over h^q(u)}{\rm d}u<+\infty$
.\vspace{0.6cm}



\end{document}